\documentclass[reqno,11pt]{article}
\usepackage{amssymb}
\usepackage{amsmath}
\usepackage{amsthm}
\usepackage{geometry} 

\newtheorem{Theorem}{Theorem } [section]
\newtheorem{lemma}[Theorem]{Lemma}

\newtheorem{remark}{Remark}

\numberwithin{equation}{section}

\newcommand{\vertiii}[1]{{\vert\kern-0.25ex\vert\kern-0.25ex\vert #1 
    \vert\kern-0.25ex\vert\kern-0.25ex\vert}}

\DeclareMathOperator{\supp}{supp}

\date{}  

\title{Dispersive blow-up for a coupled Schrödinger-fifth order KdV system}
\author{Eddye Bustamante, Jos\'e Jim\'enez Urrea and Jorge Mej\'{\i}a}

\author{
    Eddye Bustamante \\
    \textit{Universidad Nacional de Colombia, Sede Medell\'in} \\
    \texttt{eabusta0@unal.edu.co} 
    \and
    José Jiménez Urrea \\
    \textit{Universidad Nacional de Colombia, Sede Medell\'in} \\
    \texttt{jmjimene@unal.edu.co} 
    \and
    Jorge Mej\'ia \\
    \textit{Universidad Nacional de Colombia, Sede Medell\'in} \\
    \texttt{jemejia@unal.edu.co} 
}

\begin{document}

%
%

\maketitle

\begin{abstract} In this work we establish a dispersive blow-up result for the initial value problem (IVP) for the coupled Schrödinger-fifth order Korteweg-de Vries system
\begin{align*}
\left. \begin{array}{rl}
i u_t+\partial_x^2 u &\hspace{-2mm}=\alpha uv + \gamma |u|^2 u, \quad x\in\mathbb R,\quad t\in\mathbb R,\\
\partial_t v + \partial_x^5 v + \partial_x v^2&\hspace{-2mm}=\epsilon \partial_x |u|^2, \quad x\in\mathbb R,\quad t\in\mathbb R,\\
u(x,0)&\hspace{-2mm}= u_0(x), \quad v(x,0)=v_0(x).
\end{array} \right\}
\end{align*}
To achieve this, we prove a local well-posedness result in Bourgain spaces of the type $X^{s+\beta,b}\times Y^{s,b}$, along with a regularity property for the nonlinear part of the IVP solutions. This property enables the construction of initial data that leads to the dispersive blow-up phenomenon.
\end{abstract}

\setcounter{page}{1}
\pagenumbering{arabic}

\section{Introduction}

\footnotetext[1]{\textbf{2020 Mathematics Subject Classification:} 35Q53}

In this article we consider the initial value problem (IVP) for the coupled Schrödinger-fifth order Korteweg-de Vries system on the line:
\begin{align}
\left. \begin{array}{rl}
i u_t+\partial_x^2 u &\hspace{-2mm}=\alpha uv + \gamma |u|^2 u, \quad x\in\mathbb R,\quad t\in\mathbb R,\\
\partial_t v + \partial_x^5 v + \partial_x v^2&\hspace{-2mm}=\epsilon \partial_x |u|^2, \quad x\in\mathbb R,\quad t\in\mathbb R,\\
u(x,0)&\hspace{-2mm}= u_0(x), \quad v(x,0)=v_0(x),
\end{array} \right\}\label{maineq}
\end{align}

where $u=u(x,t)\in \mathbb C$, $v=v(x,t)\in\mathbb R$, and $\alpha,\gamma,\epsilon$ are real constants.\\

The system \eqref{maineq} is a particular case of the Schrödinger-Kawahara system
\begin{align}
\left. \begin{array}{rl}
i u_t+\partial_x^2 u &\hspace{-2mm}=\alpha uv + \gamma |u|^2 u,\\
\partial_t v + \beta \partial_x^3 v-\delta \partial_x^5 v + v \partial_x v&\hspace{-2mm}=\epsilon \partial_x |u|^2,\\
u(x,0)&\hspace{-2mm}= u_0(x), \quad v(x,0)=v_0(x),
\end{array} \right\}\label{maineq2}
\end{align}

where $\alpha,\beta,\gamma,\delta$ are real constants.\\

The system \eqref{maineq2} was studied by Shin in \cite{S2023}. It describes phenomena related to the interaction between short $u$ and long $v$ waves. Shin proves in this article that the IVP \eqref{maineq2} is locally well-posed when the initial data $u_0$ and $v_0$ are taken in the Sobolev spaces $H^s(\mathbb R)$ and $H^k(\mathbb R)$, respectively, with $s\geq 0$, $k\geq -2$ and $s-4\leq k<\min\{8s+1,s+2\}$, provided $\delta\neq 0$. This result is achieved by combining technics developed by Guo and Wan in \cite{GW2010} and Kato in \cite{K2011} for the Schrödinger-Korteweg-de Vries system. Furthermore, Shin establishes the global well-posedness of the IVP \eqref{maineq2} when $(u_0,v_0)\in L^2(\mathbb R)\times L^2(\mathbb R)$, under the condition $\delta\neq 0$, by adapting the method employed by Colliander, Holmer and Tzirakis in \cite{CHT2008} for the Zakharov and Klein-Gordon-Schrödinger systems in low regularity settings.\\

Although in our work we consider initial data that meet the hypotheses of Shin's article, we again establish the local well-posedness of the IVP \eqref{maineq} in the context of specific Bourgain spaces, which allow us, through their norms, to measure the regularity gains associated with the nonlinear terms of the system.\\

Let us remember that, formally, a pair $(u,v)$ is a solution of the IVP \eqref{maineq} if and only if
\begin{align}
u(t) &= U(t)u_0 -i \int_0^t U(t-t') [\alpha u(t')v(t')+\gamma |u(t')|^2 u(t') ] dt', \label{inteq1}\\
v(t) &= V(t)u_0 + \int_0^t V(t-t') [\epsilon \partial_x |u(t')|^2 - \partial_x v(t')^2] dt',
\label{inteq2}
\end{align}

where $\{U(t)\}_{t\in\mathbb R}$ and $\{V(t)\}_{t\in\mathbb R}$ are the groups of unitary operators associated to the linear partss of the Schrödinger and fifth order KdV equations, respectively, defined, through the spatial Fourier transform by 
\begin{align}
[U(t)u_0]^\wedge(\xi) &= e^{-it\xi^2} \widehat  u_0(\xi). \label{group1}\\
[V(t)v_0]^\wedge(\xi) &= e^{-it\xi^5} \widehat  v_0(\xi). \label{group2}
\end{align}

In \eqref{inteq1} and \eqref{inteq2}, $u(t)\equiv u(\cdot_x,t)$, $v(t)\equiv v(\cdot_x,t)$.\\

We are going to work in the context of the Bourgain spaces of tempered distributions in $\mathbb R_{xt}^2$, defined by
\begin{align}
X^{s,b}&\equiv X^{s,b}(\mathbb R^2):=\{u\in S'(\mathbb R^2): \|\langle \xi \rangle^s \langle \tau + \xi^2 \rangle^b \widehat  u (\xi,\tau)\|_{L^2_{\xi\tau}}<\infty\},\label{norm1}\\
Y^{s,b}&\equiv Y^{s,b}(\mathbb R^2):=\{v\in S'(\mathbb R): \|\langle \xi \rangle^s \langle \tau + \xi^5 \rangle^b \widehat  v (\xi,\tau)\|_{L^2_{\xi\tau}}<\infty\},\label{norm2}
\end{align}

where $s\in \mathbb R$, and $b=\frac12^-$, i.e. there exists $\epsilon>0$ such that $\frac12-\epsilon<b<\frac12$. Here, $\langle\cdot\rangle$ is an abbreviation for $1+|\cdot|$, $\widehat u$ and $\widehat v$ denote the Fourier transform of the distributions $u$ and $v$, respectively, and $\xi$ and $\tau$ are the variables in the frecuency space corresponding to $x$ and $t$, respectively.\\

We prove that, for suitable values of $s$, $\beta$ and $b$, and for $u_0\in H^{s+\beta}(\mathbb R)$ and $v_0\in H^s(\mathbb R)$, there exists $(u,v)\in (X^{s+\beta,b}\cap C(\mathbb R_t;H^{s+\beta}(\mathbb R_x)))\times(Y^{s,b}\cap C(\mathbb R_t;H^{s}(\mathbb R_x)))$, such that, for some $T>0$ and $t\in[0,T]$, equations \eqref{inteq1} and \eqref{inteq2} hold. Here $C(\mathbb R_t;H^{s+\beta}(\mathbb R_x))$ and $C(\mathbb R_t;H^{s}(\mathbb R_x))$ denote the spaces of continuous functions from $\mathbb R_t$ with values in $H^{s+\beta}(\mathbb R_x)$ and $H^s(\mathbb R_x)$, respectively. To achieve our result, we introduce a bump function $\eta$ of the time variable $t$ such that $\eta(\cdot_t)\in C_0^\infty(\mathbb R_t)$, $\supp \eta\subset [-1,1]$, and $\eta\equiv 1$ in $[-\frac12,\frac12]$, and we consider the integral equations

\begin{align}
u(t) &= \eta(t) U(t)u_0 -i \eta(t) \int_0^t U(t-t') F_T(u(t'),v(t')) dt', \label{inteqtr1}\\
v(t) &= \eta(t)V(t)v_0 + \eta(t)\int_0^t V(t-t') G_T(u(t'),v(t')) dt',\label{inteqtr2}
\end{align}

where
\begin{align}
F_T(u(t'),v(t')) &:= \eta(\tfrac{t'}{2T}) [\alpha u(t') v(t')+\gamma |u(t')|^2 u(t')],\label{F_T}\\ 
G_T(u(t'),v(t')) &:= \eta(\tfrac{t'}{2T}) [\epsilon \partial_x |u(t')|^2 -\partial_x v(t')^2]\label{G_T}.
\end{align}

In a precise manner, we prove the following result.
\begin{Theorem}\label{localsol} (Local well posedness of the IVP \eqref{maineq}). Let $s\geq 0$, $\frac9{20}<b<\frac12$, $\frac12<\beta<\frac52b-\frac58$ and let  $u_0\in H^{s+\beta}(\mathbb R_x)$ and $v_0\in H^s(\mathbb R_x)$. There exist $T\in(0,\frac12]$ and a unique
$$(u,v)\in (X^{s+\beta,b}\cap C(\mathbb R_t; H^{s+\beta}(\mathbb R_x)))\times (Y^{s,b}\cap C(\mathbb R_t;H^s(\mathbb R_x)))$$
solution of the integral equations \eqref{inteqtr1} and \eqref{inteqtr2}.\\

(The restriction of $u$ and $v$ to $\mathbb R\times[0,T]$ is a local solution in time of the IVP \eqref{maineq}).
\end{Theorem}

We would like to point out that the result obtained in Theorem \ref{localsol} does not include negative indices. However, the functional framework in which it is set is suitable for our purposes in this article. In contrast, Shin's result in \cite{S2023} covers a much broader range of low-regularity spaces, for which it is necessary to use Besov-type spaces.\\

In this work, our main focus is on studying dispersive blow-up for the coupled system \eqref{maineq}. From a physical perspective, this phenomenon reflects the formation of singularities caused by the interaction between short and long waves. For coupled systems of dispersive-type equations, the first result in this direction was obtained by Linares and Palacios in \cite{LP2019}, where they consider a nonlinear Schrödinger-Korteweg de Vries system and prove that the nonlinear part of the solution of the system is smoother than the solution. In this manner the existence of dispersive blow-up for the system is provided by the linear dispersive solution. Our next result is also to prove that the nonlinear part of the solution of the IVP \eqref{maineq} is smoother than the solution.\\

Some similar results have been obtained for nonlinear dispersive equations such as the KdV equation (see \cite{BS1993}), the nonlinear Schrödinger equation (see \cite{BPSS2014}) and the Zakharov-Kuznetsov equation (see \cite{LPD2021}). For related results, see also \cite{MP2023}.\\

In \cite{LP2019}, in order to prove this result, it is required, in addition to the smoothing properties of the groups associated with the linear Schrödinger and linear KdV equations, the use of the persistence property of the solution of the system in weighted Sobolev spaces.\\

For the system considered in our work, the use of Bourgain spaces and the regularity gain of the nonlinear terms of the system, measured through its norms in these spaces, are sufficient to prove that the nonlinear part of the solution is smoother than the solution. This approach, inspired by the work of Erdogan and Tzirakis \cite{ET2016}, avoids the use of the persistence property of the solution of the system in weighted Sobolev spaces, and we think that it simplifies the reasoning significantly. The precise formulation of this result is as follows.
\begin{Theorem}\label{more_reg} (Regularity gain of the nonlinear part of the solution). Let $s\geq 0$, $\frac9{20}<b<\frac12$ and $\frac12<\beta<\frac52b-\frac58$. Let $u_0\in H^{s+\beta}(\mathbb R_x)$ and $v_0\in H^s(\mathbb R_x)$. Let us take $\overline a$ and $a$ such that
$$0\leq \overline a\leq \frac52 b-\frac58-\beta \text{ and }0\leq a\leq 5\beta-\frac94.$$\\
On the other hand, let $T\in (0,\frac 12]$ and
$$(u,v)\in (X^{s+\beta,b}\cap C(\mathbb R_t;H^{s+\beta}(\mathbb R_x)))\times(Y^{s,b}\cap C(\mathbb R_t;H^s(\mathbb R_x))),$$
as in Theorem \ref{localsol}. Then
\begin{align}
\left\| \int_0^{\cdot_t} U(\cdot_t{\scriptstyle -}t')(\alpha u(t')v(t') +\gamma |u(t')|^2 u(t')) dt'\right\|_{C([0,T];H^{s+\beta+\overline a}(\mathbb R_x))} \leq C (\|u\|_{X^{s+\beta,b}} \|v\|_{Y^{s,b}} + \|u\|^3_{X^{s+\beta,b}}),\label{more_reg_eq1}
\end{align}
and
\begin{align}
\left\| \int_0^{\cdot_t} V(\cdot_t{\scriptstyle -}t')(\epsilon \partial_x|u(t')|^2 - \partial_x(v(t')^2) dt'\right\|_{C([0,T];H^{s+a}(\mathbb R_x))} \leq C (\|u\|^2_{X^{s+\beta,b}} + \|v\|^2_{Y^{s,b}}).\label{more_reg_eq2}
\end{align}
\end{Theorem}

Our main goal in this paper is to establish a dispersive blow-up result for the IVP \eqref{maineq}. In \cite{BS1993}, Bona and Saut prove the existence of dispersive blow-up for the generalized KdV equation. The idea behind the Bona and Saut proof is to observe that the nonlinear part of the solution is smoother than the solution. In this manner it is sufficient to choose initial data that gives rise to the formation of singularities for the solution of the linearized IVP. Following this idea, in \cite{LPS2017} and \cite{LP2019}, the authors establish the dispersive blow-up for the generalized KdV equation and the Schrödinger-KdV system, respectively.\\

Using the regularity gain in Theorem \ref{more_reg} and by imitating the procedure in \cite{LP2019} for the construction of the initial data, we prove the following result on dispersive blow-up for the IVP \eqref{maineq}.
\begin{Theorem}\label{blowup} (Dispersive blow-up). There exist initial data $u_0\in C^\infty(\mathbb R)\cap H^{2^-}(\mathbb R)\cap L^\infty(\mathbb R)$ and $v_0\in C^{\infty}(\mathbb R)\cap H^{\frac32^-}(\mathbb R)\cap L^{\infty}(\mathbb R)$, and $t^*\in (0,T]$ ($[0,T]$ is the time interval in which the local solution $(u,v)$ of the IVP \eqref{maineq} with initial data $u_0$ and $v_0$, is defined) such that
\begin{itemize}
\item[(i)] $u(t^*)\notin C^{1,\frac12+\epsilon}(\mathbb R)$ for some $\epsilon\in(0,\frac12)$ and
\item[(ii)] $v(t^*)\notin C^1(\mathbb R)$.
\end{itemize}
(Let us recall that $H^{s^-}(\mathbb R)=\bigcap_{s'<s} H^{s'}(\mathbb R)$. The definition of the space $C^{1,\alpha}(\mathbb R)$ is given in Section 5).
\end{Theorem}

Our work is organized as follows. In Section 2 we present the linear and nonlinear estimates needed for proving the local-well posedness of IVP \eqref{maineq}. It is important to point out the gain of regularity in the nonlinear estimates. In Section 3 we prove the local-well posedness of IVP \eqref{maineq} (proof of Theorem \ref{localsol}). In Section 4 we establish the gain of regularity of the nonlinear part of the solution of IVP \eqref{maineq} (proof of Theorem \ref{more_reg}). Finally, in Section 5, we construct the initial data $u_0$ and $v_0$ that gives rise to de dispersive blow-up (proof of Theorem \ref{blowup}).\\

Throughout the paper the letter $C$ will denote diverse constants, which may change from line to line, and whose dependence on certain parameters is clearly established in all cases.

\newpage
\section{Linear and non linear estimates}

In this section we present estimates for the norms in the spaces $X^{s,b}$ and $Y^{s,b}$ of the terms on the right hand side of the integral equations \eqref{inteqtr1} and \eqref{inteqtr2}, respectively. Proofs for some estimates are provided, while others are referenced. Throughout this section the function $\eta(\cdot_t)$ is as previously defined.\\

The proofs of the linear estimates stated in Lemmas \ref{le1}, \ref{le2} and \ref{le3} follow the ideas presented in the proofs of Lemmas 2.1, 2.2 and 2.3 in \cite{BJM2025}, respectively.

\begin{lemma}\label{le1} Let $s\in\mathbb R$, $b\in\mathbb R$. Then there exists $C>0$ such that for every function $u_0\in H^s(\mathbb R_x)$ and every function $v_0\in H^s(\mathbb R_x)$,
\begin{align}
\|\eta(\cdot_t)[U (\cdot_t)u_0](\cdot_x)\|_{X^{s,b}}&\leq C\|u_0\|_{H^s},\label{le1_eq1}\\
\|\eta(\cdot_t)[V (\cdot_t)v_0](\cdot_x)\|_{Y^{s,b}}&\leq C\|v_0\|_{H^s}.\label{le1_eq2}
\end{align}
\end{lemma}




\begin{lemma}\label{le2}
For $s\in\mathbb R$, $-\frac12<b'\leq 0\leq b\leq b'+1$ and $0<T\leq 1$, there exists $C>0$, such that for every $f\in X^{s,b'}$ and every $g\in Y^{s,b'}$
\begin{align}
\left\|\eta\left(\tfrac{\cdot_t}T \right)\int_0^{\cdot_t} [U (\cdot_t{\scriptstyle -}t')f(t')](\cdot_x)dt'\right\|_{X^{s,b}}\leq C T^{1-b+b'}\|f\|_{X^{s,b'}},\label{le2_eq1}\\
\left\|\eta\left(\tfrac{\cdot_t}T \right)\int_0^{\cdot_t} [V (\cdot_t{\scriptstyle -}t')g(t')](\cdot_x)dt'\right\|_{Y^{s,b}}\leq C T^{1-b+b'}\|f\|_{Y^{s,b'}}.\label{le2_eq2}
\end{align}
\end{lemma}

\begin{lemma}\label{le3} Let $s\in\mathbb R$, $b_1,b_2$ such that $-\frac12<b_1<b_2<\frac12$, and $T\in(0,1]$. Then there exists $C>0$ such that  for every $F\in X^{s,b_2}$ and every $G\in Y^{s,b_2}$
\begin{align}
\left\|\eta\left(\tfrac{\cdot_t}T\right)F(\cdot_x,\cdot_t) \right\|_{X^{s,b_1}}\leq CT^{b_2-b_1}\|F\|_{X^{s,b_2}},\label{le3_eq1}\\
\left\|\eta\left(\tfrac{\cdot_t}T\right)G(\cdot_x,\cdot_t) \right\|_{X^{s,b_1}}\leq CT^{b_2-b_1}\|G\|_{X^{s,b_2}}.\label{le3_eq2}
\end{align}
\end{lemma}

%
%
%
%
%

Now we present some bilinear and trilinear estimates. For a proof of Lemma \ref{ble1} we refer to Proposition 3.5 in \cite{ET2016} and Proposition 1 in \cite{ET2013}. The proof of Lemma \ref{ble2} can be found in \cite{BJM2025} (see Lemma 2.5).

\begin{lemma}\label{ble1} Let $s>\frac12$, $\frac5{12}<b<\frac12$ and $0\leq a\leq 6b-\frac52$. There exists $C>0$ such that for every $u\in X^{s,b}$, $v\in X^{s,b}$ and $w\in X^{s,b}$,
\begin{align}
\|u \overline v w\|_{X^{s+a,-b}}\leq C\|u\|_{X^{s,b}} \|v\|_{X^{s,b}}\|w\|_{X^{s,b}}.\label{ble1_eq1}
\end{align}
\end{lemma}

\begin{lemma}\label{ble2} Let $s\geq 0$,  $\frac 25< b<\frac12$ and $0\leq a\leq 10b-4$. There exists $C>0$ such that for every $v\in Y^{s,b}$ and every $w\in Y^{s,b}$
\begin{align}
\|\partial_x(vw)\|_{Y^{s+a,-b}}\leq C\|v\|_{Y^{s,b}}\|w\|_{Y^{s,b}}.\label{ble2_eq1}
\end{align}
\end{lemma}

The following calculus inequality will be used in the proof of Lemma \ref{ble3}.
\begin{lemma}\label{CI} For $\beta\geq\gamma\geq0$, and $\quad \beta+\gamma>1$ it follows that
\begin{align}
\int_{\mathbb R}\frac1{\langle x-a_1\rangle^\beta\langle x-a_2\rangle^\gamma}dx\leq C\frac{\phi_\beta(a_1-a_2)}{\langle a_1-a_2\rangle^\gamma},\label{CI_eq1}
\end{align}

where
$$\phi_\beta(a)\sim \left\{ \begin{array}{lr}
1&\text{for }\beta>1,\\
\log(1+\langle a\rangle)&\text{for }\beta=1,\\
\langle a\rangle^{1-\beta}&\text{for }\beta<1.
\end{array}\right.
$$
(For a proof of \eqref{CI_eq1}, see \cite{ET2013}).
\end{lemma}

\begin{lemma}\label{ble3} Let $s\geq 0$, $\frac9{20}<b<\frac12$ and let $a$ and $\beta$ be such that $\frac12<\beta<a\leq \frac 52b-\frac58$. There exists $C>0$ such that for every $u\in X^{s+\beta,b}$ and every $v\in Y^{s,b}$
\begin{align}
\|uv\|_{X^{s+a,-b}} \leq C \|u\|_{X^{s+\beta,b}} \|v\|_{Y^{s,b}}. \label{ble3_eq1}
\end{align}
\end{lemma}

\begin{proof} Let us observe that
$$(uv)^\wedge(\xi,\tau)=C \int_{\mathbb R_2} \widehat u(\xi_1,\tau_1)\widehat v(\xi-\xi_1,\tau-\tau_1)d\xi_1 d\tau_1.$$

Hence
\begin{align}
\notag &\|uv\|^2_{X^{s+a,-b}}\\
\notag&=C\int_{\mathbb R^2_{\xi\tau}}\langle\xi\rangle^{2(s+a)}\langle\tau+\xi^2\rangle^{-2b} \left|\int_{\mathbb R^2_{\xi_1\tau_1}}\widehat u(\xi_1,\tau_1)\widehat v(\xi-\xi_1,\tau-\tau_1)d\xi_1 d\tau_1 \right|^2d\xi d\tau\\
&=C\int_{\mathbb R^2_{\xi\tau}}\left| \int_{\mathbb R^2_{\xi_1\tau_1}} \langle\xi\rangle^{(s+a)}\langle\tau+\xi^2\rangle^{-b} \widehat u(\xi_1,\tau_1)\widehat v(\xi-\xi_1,\tau-\tau_1)d\xi_1 d\tau_1 \right|^2d\xi d\tau.\label{ble1_eq2}
\end{align}

Let $h\in L^2(\mathbb R^2_{\xi\tau})$ an arbitrary function. If we manage to prove that
\begin{align}
\notag\Bigg|\int_{\mathbb R^2_{\xi\tau}} \Bigg[\int_{\mathbb R^2_{\xi_1\tau_1}}\langle\xi\rangle^{(s+a)}\langle\tau+\xi^2\rangle^{-b} \widehat u(\xi_1,\tau_1)\widehat v(\xi-\xi_1,\tau-\tau_1) & d\xi_1 d\tau_1 \Bigg]h(\xi,\tau) d\xi d\tau\Bigg|\\
&\leq C\|u\|_{X^{s+\beta,b}}\|v\|_{Y^{s,b}}\|h\|_{L^2(\mathbb R^2)},\label{ble1_eq3}
\end{align}

then we would have, by a duality argument, that
\begin{align}
\| uv \|_{X^{s+a,-b}}\leq C\|u\|_{X^{s+\beta,b}}\|v\|_{Y^{s,b}}.\label{ble1_eq4}
\end{align}

Taking into account that, there exists $C>0$, such that for $s\geq 0$
$$\frac{\langle\xi\rangle^s}{\langle\xi_1\rangle^s\langle\xi-\xi_1\rangle^s}\leq C,$$
then, to establish \eqref{ble1_eq3}, it is enough to prove that
\begin{align}
\notag \Bigg |\int_{\mathbb R^2_{\xi\tau}} \Bigg[\int_{\mathbb R^2_{\xi_1\tau_1}} \langle\xi\rangle^{a}\langle\xi_1\rangle^s\langle\tau+\xi^2\rangle^{-b} \langle\xi-\xi_1\rangle^s|\widehat u(\xi_1,\tau_1)| & |\widehat v(\xi-\xi_1,\tau-\tau_1)|d\xi_1 d\tau_1 \Bigg] h(\xi,\tau)d\xi d\tau\Bigg|\\
&\leq C\|u\|_{X^{s+\beta,b}}\|v\|_{Y^{s,b}}\|h\|_{L^{2}(\mathbb R^2)}.\label{ble1_eq5}
\end{align}

Since
\begin{align*}
\Bigg|\int_{\mathbb R^2_{\xi\tau}} \Bigg[\int_{\mathbb R^2_{\xi_1\tau_1}} &\langle\xi\rangle^a\langle\xi_1\rangle^s\langle\tau+\xi^2\rangle^{-b} \langle\xi-\xi_1\rangle^s|\widehat u(\xi_1,\tau_1)||\widehat v(\xi-\xi_1,\tau-\tau_1)|d\xi_1 d\tau_1 \Bigg]\\
& h(\xi,\tau)d\xi d\tau\Bigg|\\
\leq\hspace{2mm} \int_{\mathbb R^2_{\xi\tau}} \int_{\mathbb R^2_{\xi_1\tau_1}} &\frac{\langle\xi\rangle^a |h(\xi,\tau)| \langle\xi_1\rangle^s \langle\xi_1\rangle^{\beta} \langle\tau_1+\xi_1^2\rangle^{b}|\widehat u(\xi_1,\tau_1)|\langle\xi-\xi_1\rangle^{s}\langle\tau-\tau_1+(\xi-\xi_1)^5\rangle^{b}}{\langle\tau+\xi^2\rangle^b\langle\tau_1+\xi_1^2\rangle^b\langle\tau-\tau_1+(\xi-\xi_1)^5\rangle^{b}\langle\xi_1\rangle^{\beta}}\\
&|\widehat v(\xi-\xi_1,\tau-\tau_1)|d\xi_1d\tau_1 d\xi d\tau\\
\leq\int_{\mathbb R^2_{\xi\tau}}\Bigg[ \int_{\mathbb R^2_{\xi_1\tau_1}}&\frac{\langle\xi\rangle^{2a} |h(\xi,\tau)|^2}{\langle\tau+\xi^2\rangle^{2b}\langle\tau_1+\xi_1^2\rangle^{2b}\langle\tau-\tau_1+(\xi-\xi_1)^5\rangle^{2b} \langle\xi_1\rangle^{2\beta}}d\xi_1 d\tau_1\Bigg]^{\frac12}\\
\Bigg[\int_{\mathbb R^2_{\xi_1\tau_1}}&\langle\xi_1\rangle^{2s+2\beta}\langle\tau_1+\xi_1^2\rangle^{2b}|\widehat u(\xi_1,\tau_1)|^2\langle\xi-\xi_1\rangle^{2s}\langle\tau-\tau_1+(\xi-\xi_1)^5\rangle^{2b}\\
&|\widehat v(\xi-\xi_1,\tau-\tau_1)|^2d\xi_1 d\tau_1 \Bigg]^{\frac12}d\xi d\tau\\
\leq\Bigg[\int_{\mathbb R^2_{\xi\tau}}\int_{\mathbb R^2_{\xi_1\tau_1}}&\frac{\langle\xi\rangle^{2a} |h(\xi,\tau)|^2}{\langle\tau+\xi^2\rangle^{2b}\langle\tau_1+\xi_1^2\rangle^{2b}\langle\tau-\tau_1+(\xi-\xi_1)^5\rangle^{2b} \langle\xi_1\rangle^{2\beta}}d\xi_1 d\tau_1 d\xi d\tau\Bigg]^{\frac12}\\
\Bigg[\int_{\mathbb R^2_{\xi\tau}}\int_{\mathbb R^2_{\xi_1\tau_1}}&\langle\xi_1\rangle^{2s+2\beta}\langle\tau_1+\xi_1^2\rangle^{2b}|\widehat u(\xi_1,\tau_1)|^2\langle\xi-\xi_1\rangle^{2s}\langle\tau-\tau_1+(\xi-\xi_1)^5\rangle^{2b}\\
&|\widehat v(\xi-\xi_1,\tau-\tau_1)|^2d\xi_1 d\tau_1d\xi d\tau \Bigg]^{\frac12}
\end{align*}
\begin{align*}
=\Bigg[ \int_{\mathbb R^2_{\xi\tau}}&|h(\xi,\tau)|^2\left(\int_{\mathbb R^2_{\xi_1\tau_1}}\frac{\langle\xi\rangle^{2a}}{\langle\tau+\xi^2\rangle^{2b}\langle\tau_1+\xi_1^2\rangle^{2b}\langle\tau-\tau_1+(\xi-\xi_1)^5\rangle^{2b} \langle\xi_1\rangle^{2\beta}}d\xi_1 d\tau_1 \right)d\xi d\tau\Bigg]^{\frac12}\\
&\|u\|_{X^{s+\beta,b}}\|v\|_{Y^{s,b}},
\end{align*}

then, to prove \eqref{ble1_eq5}, it is enough to prove that
\begin{align}
\sup_{(\xi,\tau)\in\mathbb R^2}\int_{\mathbb R^2_{\xi_1\tau_1}}\frac{\langle\xi\rangle^{2a}}{\langle\tau+\xi^2\rangle^{2b}\langle\tau_1+\xi_1^2\rangle^{2b}\langle\tau-\tau_1+(\xi-\xi_1)^5\rangle^{2b}\langle\xi_1\rangle^{2\beta}}d\xi_1d\tau_1\leq C.\label{ble1_eq6}
\end{align}

Let us observe that
\begin{align*}
&\int_{\mathbb R_{\tau_1}}\frac{\langle\xi\rangle^{2a}}{\langle\tau+\xi^2\rangle^{2b}\langle\tau_1+\xi_1^{2}\rangle^{2b}\langle\tau-\tau_1+(\xi-\xi_1)^5\rangle^{2b} \langle \xi_1\rangle^{2\beta}}d\tau_1\\
&=\frac{\langle\xi\rangle^{2a}}{\langle\tau+\xi^2\rangle^{2b} \langle \xi_1 \rangle^{2\beta}}\int_{\mathbb R_{\tau_1}}\frac1{\langle\tau_1-(-\xi_1^2)\rangle^{2b}\langle\tau_1-(\tau+(\xi-\xi_1)^5)\rangle^{2b}}d\tau_1.
\end{align*}

Using inequality \eqref{CI_eq1}, with $\beta=\gamma=2b<1$, and $\beta+\gamma=4b>1$ $(b>\frac14)$, we conclude that
\begin{align}
\notag &\int_{\mathbb R_{\xi_1}}\int_{\mathbb R_{\tau_1}}\frac{\langle\xi\rangle^{2a} }{\langle\tau+\xi^2\rangle^{2b}\langle\tau_1+\xi_1^2\rangle^{2b}\langle\tau-\tau_1+(\xi-\xi_1)^5\rangle^{2b} \langle \xi_1 \rangle^{2\beta}}d\tau_1d\xi_1\\
&\leq C \frac{\langle\xi\rangle^{2a}}{\langle\tau+\xi^2\rangle^{2b}}  \int_{\mathbb R_{\xi_1}}\frac{1}{\langle \tau+(\xi-\xi_1)^5 + \xi_1^2\rangle^{4b-1}  \langle \xi_1 \rangle^{2\beta}}d\xi_1.\label{ble1_eq7}
\end{align}

Let us consider the region $B:=\{(\xi,\tau)\in\mathbb R^2: \tau>-\frac12\xi^2\text{ or } \tau<-\frac32 \xi^2\}$.

\begin{itemize}
\item[(i)] For $(\xi,\tau)\in B$ we have that $|\tau+\xi^2|>\frac12\xi^2$. Therefore, taking into account that $\beta>\frac12$ and that $a\leq 2b$,
\begin{align*}
\frac{\langle\xi\rangle^{2a}}{\langle\tau+\xi^2\rangle^{2b}}  \int_{\mathbb R_{\xi_1}}\frac{1}{\langle \tau+(\xi-\xi_1)^5 + \xi_1^2\rangle^{4b-1}  \langle \xi_1 \rangle^{2\beta}}d\xi_1\leq C \frac{\langle \xi \rangle^{2a}}{\langle \xi^2 \rangle^{2b}} \int_{\mathbb R_{\xi_1}} \frac{d\xi_1}{\langle \xi_1 \rangle^{2\beta}} \leq C\frac{\langle \xi \rangle^{2a}}{\langle \xi^2\rangle^{2b}}\leq C.
\end{align*}

\item[(ii)] For $\tau+\xi^2=0$, we define
\begin{align*}
A_1(\xi)&:=\{\xi_1\in\mathbb R: |\xi+\xi_1|\leq \frac12|\xi-\xi_1|^4\},\\
A_2(\xi)&:=\{\xi_1\in\mathbb R: |\xi-\xi_1|^4\leq \frac12|\xi+\xi_1|\},\\
A_3(\xi)&:=\{\xi_1\in\mathbb R: \frac12|\xi-\xi_1|^4\leq |\xi+\xi_1|\leq 2|\xi-\xi_1|^4\}.
\end{align*}

Thus
\begin{align*}
\frac{\langle\xi\rangle^{2a}}{\langle\tau+\xi^2\rangle^{2b}}  \int_{\mathbb R_{\xi_1}} &\frac{1}{\langle \tau+(\xi-\xi_1)^5 + \xi_1^2\rangle^{4b-1}  \langle \xi_1 \rangle^{2\beta}}d\xi_1\\
=&\langle\xi\rangle^{2a}  \int_{A_1(\xi)}\frac{1}{\langle -\xi^2+(\xi-\xi_1)^5 + \xi_1^2\rangle^{4b-1}  \langle \xi_1 \rangle^{2\beta}}d\xi_1\\
&+\langle\xi\rangle^{2a}  \int_{A_2(\xi)}\frac{1}{\langle -\xi^2+(\xi-\xi_1)^5 + \xi_1^2\rangle^{4b-1}  \langle \xi_1 \rangle^{2\beta}}d\xi_1\\
&+\langle\xi\rangle^{2a}  \int_{A_3(\xi)}\frac{1}{\langle -\xi^2+(\xi-\xi_1)^5 + \xi_1^2\rangle^{4b-1}  \langle \xi_1 \rangle^{2\beta}}d\xi_1\\
\equiv &I + II +III.
\end{align*}

It is clear that $I+II+III$ is bounded for $|\xi|<4$. Let us assume then that $|\xi|\geq 4$, and without loss of generality, that $\xi>0$.\\

We observe that $A_1(\xi)=(-\infty,\overline\xi_0]\cup [\overline{\overline \xi}_0,+\infty)$, where $\overline \xi_0\in[0,\xi]$ and $\overline{\overline \xi}_0\in[\xi,+\infty)$ are such that
$$|\xi+\overline \xi_0|=\frac12|\xi-\overline\xi_0|^4,\quad |\xi+\overline{\overline\xi}_0|=\frac12|\xi-\overline{\overline\xi}_0|^4.$$

For $\xi_1\in A_1(\xi)$ we have that
\begin{align*}
|-\xi^2+(\xi-\xi_1)^5+\xi_1^2|&=|\xi-\xi_1||(\xi-\xi_1)^4-(\xi+\xi_1)|\geq |\xi-\xi_1|(|\xi-\xi_1|^4-|\xi+\xi_1|)\\
&\geq|\xi-\xi_1|(|\xi-\xi_1|^4-\frac12|\xi-\xi_1|^4)\geq \frac12 |\xi-\xi_1|^5.
\end{align*}

Hence, considering the change of variables $v:=\xi_1-\xi$, we obtain that
\begin{align*}
I&\leq C\langle \xi\rangle^{2a}\left[\int_{-\infty}^{\overline \xi_0} \frac{d\xi_1}{(1+|\xi-\xi_1|^5)^{4b-1}\langle \xi_1 \rangle^{2\beta}} + \int_{\overline{\overline \xi}_0}^{\infty} \frac{d\xi_1}{(1+|\xi-\xi_1|^5)^{4b-1}\langle \xi_1 \rangle^{2\beta}}\right]\\
&= C \langle \xi \rangle^{2a} \left[\int_{-\infty}^{\overline \xi_0-\xi} \frac{dv}{(1+|v|^5)^{4b-1}\langle v+\xi\rangle^{2\beta}} + \int_{\overline{\overline \xi}_0-\xi}^{\infty} \frac{dv}{(1+|v|^5)^{4b-1}\langle v+\xi\rangle^{2\beta}} \right]\\
&\leq \frac{C\langle \xi\rangle^{2a}}{(1+|\overline\xi_0-\xi|^5)^{4b-1}}+\frac{C \langle \xi \rangle^{2a}}{(1+|\overline{\overline \xi}_0-\xi|^5)^{4b-1}}
\end{align*}

Let us note that $|\xi-\overline \xi_0|^5=(|\xi-\overline \xi_0|^4)^{5/4}=(2|\xi+\overline \xi_0|)^{5/4}\geq 2^{5/4} \xi^{5/4}$ and in a similar way, that $|\xi-\overline{\overline \xi}_0|^5\geq 2^{5/4}\xi^{5/4}$.

Thereby
\begin{align*}
I\leq \frac{C\langle \xi\rangle^{2a}}{\langle \xi \rangle ^{5/4(4b-1)}}\leq C.
\end{align*}

Concerning $A_2(\xi)$, we observe that $A_2(\xi)=[\xi_0,\tilde \xi_0]$, where $\xi_0\in[0,\xi]$ and $\tilde \xi_0\in[\xi,+\infty)$ are such that
$$\frac12|\xi+\xi_0|=|\xi-\xi_0|^4,\quad \frac12|\xi+\tilde \xi_0|=|\xi-\tilde \xi_0|^4.$$

For $\xi_1\in A_2(\xi)$ we have that
\begin{align*}
|(\xi+\xi_1)-(\xi-\xi_1)^4|& \geq |\xi+\xi_1|-\frac12 |\xi+\xi_1|=\frac12|\xi+\xi_1|.
\end{align*}

Therefore, taking into account that $\xi-\xi^{1/4}<\xi_0<\xi<\tilde\xi_0<\xi+2\xi^{1/4}$ for $\xi>4$, we obtain that
\begin{align*}
II&=\langle \xi\rangle^{2a}\int_{\xi_0}^{\tilde\xi_0} \frac{d\xi_1}{\langle |\xi-\xi_1||(\xi+\xi_1)-(\xi-\xi_1)^4| \rangle^{4b-1} \langle \xi_1 \rangle^{2\beta}}\\
&\leq C\langle \xi\rangle^{2a} \int_{\xi_0}^{\tilde \xi_0} \frac{d\xi_1}{\langle |\xi-\xi_1| |\xi+\xi_1| \rangle^{4b-1} \langle \xi_1\rangle^{2\beta}}\\
&\leq C\langle \xi\rangle^{2a} \int_{\xi-\xi^{1/4}}^{\xi} \frac{d\xi_1}{\langle |\xi-\xi_1| |\xi+\xi_1| \rangle^{4b-1} \langle \xi_1\rangle^{2\beta}}+C\langle \xi\rangle^{2a} \int_{\xi}^{\xi+2\xi^{1/4}} \frac{d\xi_1}{\langle |\xi-\xi_1| |\xi+\xi_1| \rangle^{4b-1} \langle \xi_1\rangle^{2\beta}}\\
&=IIa+IIb.
\end{align*}

We only provide an estimate for $IIa$, since the estimation of $IIb$ is similar. Considering the change of variables $v=\xi-\xi_1$, and $u=1+2\xi v -v^2$, and taking into account that $\langle \xi_1 \rangle\sim \langle \xi \rangle$, and that $\frac\xi2<\xi-v<\xi$, we have
\begin{align*}
IIa&\leq C\langle \xi\rangle^{2a} \int_0^{\xi^{1/4}} \frac{dv}{[1+v(2\xi-v)]^{4b-1}\langle \xi_1\rangle^{2\beta}}\leq C \frac{\langle \xi \rangle^{2a}}{\langle \xi \rangle^{2\beta}} \int_1^{1+2\xi^{5/4}-\xi^{1/2}} \frac{du}{u^{4b-1}(\xi-v)}\\
&\leq C \frac{\langle \xi \rangle^{2a}}{\xi \langle \xi\rangle^{2\beta}} \frac{u^{2-4b}}{2-4b} \Bigg|_{u=1}^{u=1+2\xi^{5/4}-\xi^{1/2}}\leq C\frac{\langle \xi \rangle^{2a}\xi^{5/2-5b}}{(2-4b)\xi \langle \xi \rangle^{2\beta}}\leq C \frac{\langle \xi\rangle^{2a+3/2-5b-2\beta}}{(2-4b)}\leq \frac{C}{(2-4b)}\\
&\leq C.
\end{align*}

Then, we can affirm that $II\leq C$.\\

Now we consider $A_3(\xi)$. Let us observe that $A_3(\xi)=[\eta_0,\overline \eta_0]\cup[\delta_0,\overline \delta_0]$, where $\eta_0,\overline \eta_0\in[0,\xi]$, and $\delta_0,\overline \delta_0\in[\xi,+\infty)$ are such that
$$|\xi+\eta_0|=\frac12|\xi-\eta_0|^4,\quad |\xi+\overline \eta_0|=2|\xi-\overline \eta_0|^4,\quad |\xi+\delta_0|=2|\xi-\delta_0|^4,\quad |\xi+\overline\delta_0|=\frac12 |\xi-\overline \delta_0|^4.$$

Let us choose $\xi_0\in[\eta_0,\overline \eta_0]$ and $\tilde \xi_0\in[\delta_0,\overline \delta_0]$ such that
$$|\xi+\xi_0|=|\xi-\xi_0|^4,\quad |\xi+\tilde \xi_0|=|\xi-\tilde\xi_0|^4.$$
Therefore
$$A_3(\xi)=[\eta_0,\xi_0]\cup[\xi_0,\overline\eta_0]\cup[\delta_0,\tilde\xi_0]\cup[\tilde\delta_0,\overline\delta_0].$$

We will estimate only
$$\langle\xi\rangle^{2a}  \int_{[\eta_0,\xi_0]}\frac{1}{\langle -\xi^2+(\xi-\xi_1)^5 + \xi_1^2\rangle^{4b-1}  \langle \xi_1 \rangle^{2\beta}}d\xi_1,$$
since the estimation of the integrals over $[\xi_0,\overline\eta_0]$, $[\delta_0,\tilde\xi_0]$ and $[\tilde\delta_0,\overline\delta_0]$ is analogous.\\

Taking into account that $\langle \xi_1\rangle \sim \langle \xi \rangle$, that $\xi-\xi_1\geq \frac1{2^{1/4}} \xi^{1/4}$ for $\xi_1\in[\eta_0,\xi_0]$, and considering the change of variables $v=(\xi-\xi_1)^4-(\xi-\xi_1)$, it follows that
\begin{align*}
\langle\xi\rangle^{2a}&  \int_{[\eta_0,\xi_0]}\frac{1}{\langle -\xi^2+(\xi-\xi_1)^5 + \xi_1^2\rangle^{4b-1}  \langle \xi_1 \rangle^{2\beta}}d\xi_1\\
 &\leq C\frac{\langle \xi \rangle^{2a}}{\langle \xi \rangle^{2\beta}} \int_{\eta_0}^{\xi_0} \frac{d\xi_1}{(1+(\xi-\xi_1)[(\xi-\xi_1)^4-(\xi+\xi_1)])^{4b-1}}\\
 & \leq C\frac{\langle \xi \rangle^{2a}}{\langle \xi \rangle^{2\beta}} \int_{\eta_0}^{\xi_0} \frac{d\xi_1}{(1+\frac1{2^{1/4}}\xi^{1/4} [(\xi-\xi_1)^4-(\xi+\xi_1)])^{4b-1}}\\
 & \leq C\frac{\langle \xi \rangle^{2a}}{\langle \xi \rangle^{2\beta}} \int_{0}^{|\xi+\eta_0|} \frac{dv}{(1+\frac1{2^{1/4}}\xi^{1/4}v)^{4b-1}(4(\xi-\xi^1)^3+1)}\\
 &\leq \frac{C}{2-4b}\frac{\langle \xi\rangle^{2a}}{\langle \xi \rangle^{2\beta}} \frac1{\xi^{3/4}} \frac{(1+\frac{\xi^{1/4}}{2^{1/4}})^{2-4b}}{\xi^{1/4}}\Bigg|_{v=0}^{v=|\xi+\eta_0|}\leq \frac C{2-4b} \xi^{2a-2\beta-1 }(\xi^{1/4}\xi)^{2-4b}\leq \frac{C}{2-4b},
\end{align*}

if $2a-2\beta-1+\frac52-5b\leq 0$, i.e. if $a-\beta\leq \frac52 b-\frac34$, which is true because $a-\beta\leq\frac52 b -\frac98$.\\

This allows to assure that $III\leq C$. As a result, we can conclude that, for $\tau=-\xi^2$,
\begin{align}
\notag\frac{\langle\xi\rangle^{2a}}{\langle\tau+\xi^2\rangle^{2b}}  \int_{\mathbb R_{\xi_1}} &\frac{1}{\langle \tau+(\xi-\xi_1)^5 + \xi_1^2\rangle^{4b-1}  \langle \xi_1 \rangle^{2\beta}}d\xi_1\\
&=\langle\xi\rangle^{2a} \int_{\mathbb R_{\xi_1}} \frac{1}{\langle -\xi^2+(\xi-\xi_1)^5 + \xi_1^2\rangle^{4b-1}  \langle \xi_1 \rangle^{2\beta}}d\xi_1\leq C,\label{Est_con_ig}
\end{align}

\item[(iii)] To complete the estimation for the whole plane, it is left to consider the region
$$B^c\setminus \{(\xi,\tau)\in\mathbb R^2: \tau=- \xi^2\}= \{(\xi,\tau)\in\mathbb R^2: \tau=\alpha \xi^2, \, \alpha \in[-\tfrac32,-\tfrac12],\, \alpha\neq -1\}.$$
For $(\xi,\tau)\in B^c\setminus \{(\xi,\tau)\in\mathbb R^2:\tau=-\xi^2\}$, we observe that
\begin{align*}
|\alpha \xi^2 + (\xi-\xi_1)^5+\xi_1^2| & = |\alpha\xi^2 + \xi^2-\xi^2 + (\xi-\xi_1)^5+\xi_1^2|\\
& \geq |-\xi^2 + (\xi-\xi_1)^5 + \xi_1^2| - |(\alpha+1)\xi^2|.
\end{align*}

Let us consider
$$B_1(\xi):=\{\xi_1\in\mathbb R: -|(\alpha+1)\xi^2|\geq -\frac12 |-\xi^2+(\xi-\xi_1)^5+\xi_1^2| \}.$$
For $\xi_1\in B_1(\xi)$, it is clear that
$$\langle \alpha \xi^2 + (\xi-\xi_1)^5 + \xi_1^2 \rangle \geq \frac12 \langle -\xi^2+(\xi-\xi_1)^5+\xi_1^2 \rangle.$$
Therefore, in view of the estimation made in \eqref{Est_con_ig},
\begin{align*}
\frac{\langle\xi\rangle^{2a}}{\langle\tau+\xi^2\rangle^{2b}}  \int_{B_1(\xi)} &\frac{1}{\langle \tau+(\xi-\xi_1)^5 + \xi_1^2\rangle^{4b-1}  \langle \xi_1 \rangle^{2\beta}}d\xi_1\\
&\leq \frac{ \langle\xi\rangle^{2a}} {\langle\tau+\xi^2\rangle^{2b}} \int_{B_1(\xi)} \frac{2^{4b-1}}{\langle -\xi^2 + (\xi-\xi_1)^5+\xi_1^2 \rangle^{4b-1}  \langle \xi_1 \rangle^{2\beta}}d\xi_1\leq C.
\end{align*}

On the other hand, for $\xi_1\in \mathbb R\setminus B_1(\xi)$, we have that
\begin{align*}
-|(\alpha+1)\xi^2|&<-\frac12|-\xi^2+(\xi-\xi_1)^5+\xi_1^2|,\\
|(\alpha+1)\xi^2|&>\frac12|-\xi^2+(\xi-\xi_1)^5+\xi_1^2|,\\
\langle (\alpha+1)\xi^2 \rangle&>\frac12 \langle -\xi^2+(\xi-\xi_1)^5+\xi_1^2 \rangle,\\
\langle (\alpha+1)\xi^2 \rangle^{2b}& > \langle (\alpha+1)\xi^2 \rangle^{4b-1}\geq \frac1{2^{4b-1}} \langle -\xi^2+(\xi-\xi_1)^5+\xi_1^2 \rangle^{4b-1}.
\end{align*}

Hence, using again the estimation made in \eqref{Est_con_ig},
\begin{align*}
\frac{\langle\xi\rangle^{2a}}{\langle\tau+\xi^2\rangle^{2b}}  \int_{\mathbb R\setminus B_1(\xi)} &\frac{1}{\langle \tau+(\xi-\xi_1)^5 + \xi_1^2\rangle^{4b-1}  \langle \xi_1 \rangle^{2\beta}}d\xi_1\\
&\leq  \langle\xi\rangle^{2a}  \int_{\mathbb R\setminus B_1(\xi)} \frac{1}{\langle (\alpha+1)\xi^2\rangle^{2b}  \langle \xi_1 \rangle^{2\beta}}d\xi_1\\
&\leq \langle\xi\rangle^{2a}  \int_{\mathbb R\setminus B_1(\xi)} \frac{2^{4b-1}}{\langle -\xi^2+(\xi-\xi_1)^5+\xi_1^2 \rangle^{4b-1}  \langle \xi_1 \rangle^{2\beta}}d\xi_1\leq C.
\end{align*}
Consequently, for $(\xi,\tau)\in B^c\setminus \{(\xi,\tau)\in \mathbb R^2:\tau=-\xi^2\}$, we obtain
\begin{align*}
\frac{\langle\xi\rangle^{2a}}{\langle\tau+\xi^2\rangle^{2b}}  \int_{\mathbb R_{\xi_1}} &\frac{1}{\langle \tau+(\xi-\xi_1)^5 + \xi_1^2\rangle^{4b-1}  \langle \xi_1 \rangle^{2\beta}}d\xi_1\leq C.\\
\end{align*}
\end{itemize}

From (i), (ii), (iii), and inequality \eqref{ble1_eq7}, we conclude that \eqref{ble1_eq6} holds. Thus, Lemma \ref{ble3} is proven.
\end{proof}

\begin{lemma}\label{ble4} Let $s\geq 0$, $\frac38<b<\frac12$, $\frac12<\beta\leq 2b-\frac14$ and $0\leq a \leq 5\beta-\frac94$. There exists $C>0$ such that for every $u\in X^{s+\beta,b}$ and every $v\in X^{s+\beta,b}$,
\begin{align}
\|\partial_x(u\overline v)\|_{Y^{s+a,-b}} \leq C \|u\|_{X^{s+\beta,b}} \|v\|_{X^{s+\beta,b}}. \label{ble4_eq1}
\end{align}
\end{lemma}

\begin{proof} Following the same ideas as in the proof of Lemma \ref{ble3}, we see that, to prove \eqref{ble4_eq1}, it suffices to show that
\begin{align}
\sup_{(\xi,\tau)\in\mathbb R^2}\int_{\mathbb R^2_{\xi_1\tau_1}}\frac{|\xi|^2\langle\xi\rangle^{2a}}{\langle\tau+\xi^5\rangle^{2b} \langle \tau_1+\xi_1^2 \rangle^{2b} \langle \xi-\xi_1 \rangle^{2\beta}  \langle \tau-\tau_1+(\xi-\xi_1)^2\rangle^{2b}\langle\xi_1\rangle^{2\beta}}d\xi_1d\tau_1\leq C.\label{ble4_eq2}
\end{align}

As in the proof of Lemma \ref{ble3}, using inequality \eqref{CI_eq1}, we observe that
\begin{align}
\notag &\int_{\mathbb R^2_{\xi_1\tau_1}}\frac{|\xi|^2\langle\xi\rangle^{2a}}{\langle\tau+\xi^5\rangle^{2b} \langle \tau_1+\xi_1^2 \rangle^{2b} \langle \xi-\xi_1 \rangle^{2\beta}  \langle \tau-\tau_1+(\xi-\xi_1)^2\rangle^{2b}\langle\xi_1\rangle^{2\beta}}d\xi_1d\tau_1\\
&\leq C \frac{|\xi|^2 \langle\xi\rangle^{2a}}{\langle\tau+\xi^5\rangle^{2b}}  \int_{\mathbb R_{\xi_1}}\frac{1}{ \langle \xi-\xi_1\rangle^{2\beta}\langle \tau + \xi_1^2+(\xi-\xi_1)^2\rangle^{4b-1}  \langle \xi_1 \rangle^{2\beta}}d\xi_1.\label{ble4_eq3}
\end{align}

Let us consider the region $\tilde B:=\{(\xi,\tau)\in\mathbb R^2: |\tau+\xi^5|>\frac12 |\xi|^5\}$.
\begin{itemize}
\item[(i)] For $(\xi,\tau)\in \tilde B$, taking into account that $\beta>\frac12$, that $a\leq 5\beta-\frac94$ and that $b>\frac14$, it follows that
\begin{align*}
\frac{|\xi|^2 \langle\xi\rangle^{2a}}{\langle\tau+\xi^5\rangle^{2b}} & \int_{\mathbb R_{\xi_1}}\frac{1}{ \langle \xi-\xi_1\rangle^{2\beta}\langle \tau + \xi_1^2+(\xi-\xi_1)^2\rangle^{4b-1}  \langle \xi_1 \rangle^{2\beta}}d\xi_1 \leq C\frac{|\xi|^2 \langle \xi \rangle^{2a}}{\langle \xi^5 \rangle^{2b}} \int_{\mathbb R_{\xi_1}} \frac{d\xi_1}{\langle \xi_1 \rangle^{2\beta}}\leq C.
\end{align*}

\item[(ii)] For $\tau+\xi^5=0$, we define
\begin{align*}
\tilde A_1(\xi)&:=\{\xi_1\in\mathbb R: \xi_1^2+(\xi-\xi_1)^2<\frac12 |\xi|^5\},\\
\tilde A_2(\xi)&:=\{\xi_1\in\mathbb R: |\xi|^5<\frac12[\xi_1^2+(\xi-\xi_1)^2]\},\\
\tilde A_3(\xi)&:=\{\xi_1\in\mathbb R: \xi\notin A_1(\xi)\text{ and }\xi_1\notin A_2(\xi) \}.
\end{align*}

Therefore
\begin{align*}
\frac{|\xi|^2 \langle\xi\rangle^{2a}}{\langle\tau+\xi^5\rangle^{2b}} & \int_{\mathbb R_{\xi_1}}\frac{1}{ \langle \xi-\xi_1\rangle^{2\beta}\langle \tau + \xi_1^2+(\xi-\xi_1)^2\rangle^{4b-1}  \langle \xi_1 \rangle^{2\beta}}d\xi_1\\
=&|\xi|^2 \langle\xi\rangle^{2a}   \int_{\tilde A_1(\xi)}\frac{1}{ \langle \xi-\xi_1\rangle^{2\beta}\langle -\xi^5 + \xi_1^2+(\xi-\xi_1)^2\rangle^{4b-1}  \langle \xi_1 \rangle^{2\beta}}d\xi_1\\
&+ |\xi|^2 \langle\xi\rangle^{2a}   \int_{\tilde A_2(\xi)}\frac{1}{ \langle \xi-\xi_1\rangle^{2\beta}\langle -\xi^5 + \xi_1^2+(\xi-\xi_1)^2\rangle^{4b-1}  \langle \xi_1 \rangle^{2\beta}}d\xi_1\\
&+ |\xi|^2 \langle\xi\rangle^{2a}   \int_{\tilde A_3(\xi)}\frac{1}{ \langle \xi-\xi_1\rangle^{2\beta}\langle -\xi^5 + \xi_1^2+(\xi-\xi_1)^2\rangle^{4b-1}  \langle \xi_1 \rangle^{2\beta}}d\xi_1\\
\equiv &I+II +III.
\end{align*}
It is clear that for $\xi$ in a bounded set the expression $I+II+III$ is bounded by a constant $C$. Without loss of generality, we assume that $\xi>0$ for sufficiently large values of $\xi$.\\

We observe that if $\xi_1\in \tilde A_1(\xi)$, then
\begin{align*}
|-\xi^5+\xi_1^2+(\xi-\xi_1)^2|&\geq \frac12|\xi|^5,\\
\langle -\xi^5+\xi_1^2+(\xi-\xi_1)^2\rangle &>\frac12\langle \xi^5\rangle,
\end{align*}
which implies, taking into account that $\beta>\frac12$, $a \leq 5\beta-\frac94$, and $\beta\leq 2b-\frac14$, that
\begin{align*}
I&\leq C \frac{|\xi|^2 \langle \xi\rangle^{2a}}{\langle \xi\rangle^{20b-5}}\int_{\tilde A_1(\xi)} \frac{d\xi_1}{\langle \xi_1 \rangle^{2\beta} \langle \xi-\xi_1\rangle^{2\beta}}\leq C.
\end{align*}
Let us note that $A_1(\xi)=[\xi_0,\overline \xi_0]$, where $\xi_0<0<\overline \xi_0$ are such that
$$\xi_0+(\xi-\xi_0)^2=(\overline \xi_0)^2+(\xi-\overline \xi_0)^2=\frac12|\xi|^5,$$
i.e.,
$$\xi_0=\frac\xi2-\frac12\sqrt{|\xi|^5-\xi^2},\quad \overline\xi_0=\frac\xi2+\frac12\sqrt{|\xi|^5-\xi^2}.$$
In this way, for $|\xi|$ sufficiently large, $|\xi_0|$ and $|\overline\xi_0|$ are of the same order as $|\xi|^{5/2}$.\\

Concerning $\tilde A_2(\xi)$, we observe that $\tilde A_2(\xi)=(-\infty,\tilde \xi_0)\cup(\overline{\tilde \xi}_0,+\infty)$, where $\tilde \xi_0<0<\overline{\tilde\xi}_0$ are such that
$$\tilde \xi_0^2+(\xi-\tilde \xi_0)^2=\overline{\tilde \xi}_0^2+(\xi-\overline{\tilde \xi}_0)^2=2|\xi|^5,$$
i.e.,
$$\tilde \xi_0=\frac\xi2-\frac12\sqrt{4|\xi|^5-\xi^2},\quad \overline{\tilde \xi}_0=\frac\xi2+\frac12\sqrt{4|\xi|^5-\xi^2}.$$
We conclude that, for $|\xi|$ sufficiently large, $|\tilde\xi_0|$ and $|\overline{\tilde\xi}_0|$ are of the same order as $|\xi|^{5/2}$.\\

For $\xi_1\in \tilde A_2(\xi)$,
$$\langle -\xi^5+\xi_1^2+(\xi-\xi_1)^2\rangle >\frac12\langle\xi_1^2+(\xi-\xi_1)^2\rangle,$$
hence
\begin{align*}
II\leq &C|\xi|^2\langle \xi\rangle^{2a} \int_{\tilde A_2(\xi)} \frac{d\xi_1}{ \langle \xi-\xi_1\rangle^{2\beta}\langle -\xi^5 + \xi_1^2+(\xi-\xi_1)^2\rangle^{4b-1}  \langle \xi_1 \rangle^{2\beta}}\\
\leq& C|\xi|^2\langle \xi\rangle^{2a} \int_{\tilde A_2(\xi)} \frac{d\xi_1}{ \langle \xi-\xi_1\rangle^{2\beta}\langle \xi_1^2+(\xi-\xi_1)^2\rangle^{4b-1}  \langle \xi_1 \rangle^{2\beta}}\\
\leq& C|\xi|^2\langle \xi\rangle^{2a} \int_{\tilde A_2(\xi)} \frac{d\xi_1}{ \langle \xi-\xi_1\rangle^{4\beta+8b-2}}\\
=&C|\xi|^2\langle \xi\rangle^{2a} \int_{\overline{\tilde \xi}_0}^{+\infty} \frac{d\xi_1}{ \langle \xi-\xi_1\rangle^{4\beta+8b-2}}+C|\xi|^2\langle \xi\rangle^{2a} \int_{-\infty}^{\tilde\xi_0} \frac{d\xi_1} { \langle \xi-\xi_1\rangle^{4\beta+8b-2}}\equiv IIa + IIb.
\end{align*}
We only provide an estimate for $IIa$, since the estimation of $IIb$ is similar. Considering the change of variables $u=\xi-\xi_1$, and taking into account that $a\leq 5\beta-\frac{9}4$, and $\frac 38<b$, we have

\begin{align*}
IIa &= |\xi|^2\langle \xi \rangle^{2a} \int_{\overline{\tilde \xi}_0-\xi}^{+\infty} \frac {du}{(1+u)^{4\beta+8b-2}} \leq |\xi|^2\langle \xi \rangle^{2a} \int_{C|\xi|^{5/2}}^{+\infty} \frac {du}{(1+u)^{4\beta+8b-2}}\\
 &=|\xi|^2\langle \xi \rangle^{2a}\frac1{3-4\beta-8b} (1+u)^{3-4\beta-8b} \Bigg|_{u=C|\xi|^{5/2}}^{u=+\infty}= C\frac1{4\beta+8b-3} \frac{ |\xi|^2\langle \xi \rangle^{2a}}{\langle |\xi|^{5/2} \rangle^{4\beta+8b-3}}\leq C.
\end{align*}

Then, we can affirm that $II\leq C$.\\

Now we consider $\tilde A_3(\xi)$. Let us observe that
$$\tilde A_3(\xi)=\{\xi_1\in\mathbb R: \frac12 |\xi|^5<\xi_1^2+(\xi-\xi_1)^2<2|\xi|^5\}=(\tilde \xi_0,\xi_0)\cup(\overline \xi_0,\overline{\tilde \xi}_0).$$

Sinice for $\xi_1\in(\tilde\xi_0,\xi_0)\cup(\overline\xi_0,\overline{\tilde \xi}_0)$, both $|\xi_1|$ and $|\xi-\xi_1|$ are of the same order as $|\xi|^{5/2}$, and the measure of $(\tilde \xi_0,\xi_0)\cup(\overline \xi_0,\overline{\tilde \xi}_0)$ is also of the same order as $|\xi|^{5/2}$, and taking into account that $a\leq 5\beta-\frac94$, it follows that
\begin{align*}
III & \leq C \frac{|\xi|^2 \langle \xi \rangle^{2a}}{\langle\xi\rangle^{10\beta}} |\xi|^{5/2}\leq C.
\end{align*}
Therefore, we can conclude that, for $\tau=-\xi^5$,
\begin{align}
\frac{|\xi|^2 \langle\xi\rangle^{2\alpha}}{\langle\tau+\xi^5\rangle^{2b}}  \int_{\mathbb R_{\xi_1}}\frac{1}{ \langle \xi-\xi_1\rangle^{2\beta}\langle \tau + \xi_1^2+(\xi-\xi_1)^2\rangle^{4b-1}  \langle \xi_1 \rangle^{2\beta}}d\xi_1\leq C\label{ble4_eq4}
\end{align}

\item[(iii)] To complete the estimation for the whole plane, it is left to consider the region
$$\tilde B^c\setminus \{(\xi,\tau)\in\mathbb R^2: \tau=- \xi^5\}= \{(\xi,\tau)\in\mathbb R^2: \tau=\theta \xi^5, \, \theta \in[-\tfrac32,-\tfrac12],\, \theta\neq -1\}.$$

For $(\xi,\tau)\in \tilde B^c\setminus \{(\xi,\tau)\in\mathbb R^2: \tau=- \xi^5\}$, we observe that
\begin{align*}
|\theta \xi^5+\xi_1^2+(\xi-\xi_1)^2|&=|\theta\xi^5+\xi^5-\xi^5+\xi_1^2+(\xi-\xi_1)^2|\\
&\geq |-\xi^5+\xi_1^2+(\xi-\xi_1)^2|-|(1+\theta)\xi^5|.
\end{align*}

Let us consider
\begin{align*}
\tilde B_1(\xi):=\{\xi_1\in\mathbb R: - |(1+\theta)\xi^5|\geq -\frac12 |-\xi^5+\xi_1^2+(\xi-\xi_1)^2|\}.
\end{align*}
For $\xi_1\in \tilde B_1(\xi)$, it is clear that
$$\langle \theta\xi^5 + \xi_1^2+(\xi-\xi_1)^2 \rangle \geq \frac12 \langle - \xi^5+\xi_1^2+(\xi-\xi_1)^2\rangle.$$

Therefore, in view of the estimation made in \eqref{ble4_eq4},
\begin{align*}
\frac{|\xi|^2 \langle\xi\rangle^{2a}}{\langle\tau+\xi^5\rangle^{2b}} & \int_{\tilde B_1(\xi)}\frac{1}{ \langle \xi-\xi_1\rangle^{2\beta}\langle \tau + \xi_1^2+(\xi-\xi_1)^2\rangle^{4b-1}  \langle \xi_1 \rangle^{2\beta}}d\xi_1\\
&\leq 2^{4b-1} |\xi|^2 \langle \xi \rangle^{2a} \int_{\tilde B_1(\xi)} \frac{1}{ \langle \xi-\xi_1\rangle^{2\beta}\langle -\xi^5 + \xi_1^2+(\xi-\xi_1)^2\rangle^{4b-1}  \langle \xi_1 \rangle^{2\beta}}d\xi_1\leq C.
\end{align*}

On the other hand, for $\xi_1\in\mathbb R\setminus \tilde B_1(\xi)$, since $2b>4b-1$, we have that
\begin{align*}
-|(1+\theta)\xi^5|&<-\frac12|-\xi^5+(\xi-\xi_1)^2+\xi_1^2|,\\
|(1+\theta)\xi^5|&>\frac12|-\xi^5+(\xi-\xi_1)^2+\xi_1^2|,\\
\langle (1+\theta)\xi^5 \rangle&>\frac12 \langle -\xi^5+(\xi-\xi_1)^2+\xi_1^2 \rangle,\\
\langle (1+\theta)\xi^5 \rangle^{2b}& > \langle (1+\theta)\xi^5 \rangle^{4b-1}\geq \frac1{2^{4b-1}} \langle -\xi^5+(\xi-\xi_1)^2+\xi_1^2 \rangle^{4b-1}.
\end{align*}

Hence, using again the estimation made in \eqref{ble4_eq4}
\begin{align*}
\frac{|\xi|^2 \langle\xi\rangle^{2a}}{\langle\tau+\xi^5\rangle^{2b}} & \int_{\mathbb R\setminus \tilde B_1(\xi)}\frac{1}{ \langle \xi-\xi_1\rangle^{2\beta}\langle \tau + \xi_1^2+(\xi-\xi_1)^2\rangle^{4b-1}  \langle \xi_1 \rangle^{2\beta}}d\xi_1\\
&\leq 2^{4b-1} |\xi|^2 \langle \xi \rangle^{2a}  \int_{\mathbb R\setminus \tilde B_1(\xi)}\frac{1}{ \langle \xi-\xi_1\rangle^{2\beta}\langle -\xi^5 + \xi_1^2+(\xi-\xi_1)^2\rangle^{4b-1}  \langle \xi_1 \rangle^{2\beta}}d\xi_1\leq C.
\end{align*}
\end{itemize}

From (i), (ii), (iii), and inequality \eqref{ble4_eq3}, we conclude that \eqref{ble4_eq2} holds. Thus, Lemma \ref{ble4} is proven.

\end{proof}

\section{Existence of a local solution in time of IVP \eqref{maineq}. Proof of Theorem \ref{localsol}}

We begin by defining, for $u\in X^{s+\beta,b}$ and $v\in Y^{s,b}$, the operators
\begin{align*}
\Phi_T (u,v)(t) := \eta(t) U(t) u_0 - i \eta(t) \int_0^t U(t-t') F_T(u(t'),v(t')) dt',\\
\Psi_T (u,v)(t) := \eta(t) V(t) v_0 + \eta(t) \int_0^t V(t-t') G_T(u(t'),v(t')) dt'.
\end{align*}

We will consider the space of functions
$$H:=X^{s+\beta,b}\times Y^{s,b},$$
with the associated norm
$$\|(u,v)\|_H:=\|u\|_{X^{s+\beta,b}}+\|v\|_{Y^{s,b}},$$

and we will prove that there exists $T\in (0,\frac12]$ such that the operator defined by
$$\Gamma_T(u,v)(t):=(\Phi_T(u,v)(t);\Psi_T(u,v)(t))$$

has a fixed point in $H$.\\

Recall that, by inequalities \eqref{le1_eq1} and \eqref{le1_eq2} in Lemma \ref{le1}, we have that
\begin{align}
\|\eta(\cdot_t)[U (\cdot_t)u_0](\cdot_x)\|_{X^{s+\beta,b}}&\leq C\|u_0\|_{H^{s+\beta}},\label{lst_eq1}\\
\|\eta(\cdot_t)[V (\cdot_t)v_0](\cdot_x)\|_{Y^{s,b}}&\leq C\|v_0\|_{H^{s}},\label{lst_eq2}
\end{align}

and by using Lemma \ref{le2} with $T=1$, $b':=-b^*$, where $b^*$ is such that $b<b^*<\frac12$,
\begin{align}
\left\|\eta (\cdot_t)\int_0^{\cdot_t} U (\cdot_t{\scriptstyle -}t')F_T(u(t'),v(t')) dt'\right\|_{X^{s+\beta,b}}&\leq C \|F_T(u,v)\|_{X^{s+\beta,-b^*}},\label{lst_eq3}\\
\left\|\eta (\cdot_t)\int_0^{\cdot_t} V (\cdot_t{\scriptstyle -}t')G_T(u(t'),v(t')) dt'\right\|_{Y^{s,b}}&\leq C \|G_T(u,v)\|_{Y^{s,-b^*}}.\label{lst_eq4}
\end{align}

Since $-\frac12<-b^*<-b<\frac12$, by Lemma \ref{le3}, with $2T\leq 1$, and Lemmas \ref{ble1}, \ref{ble2}, \ref{ble3} and \ref{ble4}, we have that
\begin{align}
\|F_T(u,v)\|_{X^{s+\beta,-b^*}} & \leq CT^{-b+b^*} \| \alpha uv + \gamma |u|^2 u \|_{X^{s+\beta,-b}}\leq C T^{-b+b^*} ( \|u\|_{X^{s+\beta,b}} \|v\|_{Y^{s,b}}+ \|u\|^3_{X^{s+\beta,b}}),\label{lst_eq5}
\end{align}

and
\begin{align}
\|G_T(u,v)\|_{Y^{s,-b^*}} & \leq CT^{-b+b^*} \| \epsilon \partial_x |u|^2 - \partial_x v^2 \|_{Y^{s,-b}}\leq C T^{-b+b^*} ( \|u\|^2_{X^{s+\beta,b}} + \|v\|^2_{Y^{s,b}}).\label{lst_eq6}
\end{align}

Therefore, taking into account \eqref{lst_eq1} to \eqref{lst_eq6}, we can affirm that

\begin{align}
\notag\|\Gamma_T(u,v) \|_{H} = & \| \Phi_T(u,v)\|_{X^{s+\beta,b}} + \|\Psi_T(u,v)\|_{Y^{s,b}}\\
\notag\leq &C \|u_0\|_{H^{s+\beta}} +C\|v_0\|_{H^s}\\
&+C T^{-b+b^*} (\|u\|_{X^{s+\beta,b}}\|v\|_{Y^{s,b}} +\|u\|^3_{X^{s+\beta,b}} + \|u^2_{X^{s+\beta,b}}\| + \|v\|^2_{Y^{s,b}}). \label{lst_eq7}
\end{align}

Let us consider the closed ball $B_{H}(0,R)\subset H$,  centered at 0, with radius $R:=2C(\|u_0\|_{H^s} + \|v_0\|_{H^s})>0$. Hence, from \eqref{lst_eq7}, it follows that if $(u,v)\in B_{H}(0,R)$,
$$\|\Gamma_T (u,v) \|_{H} \leq \frac R2 + C T^{-b+b^*} (R^2+R^3+R^2+R^2).$$

Choosing $T>0$ such that
\begin{align}
CT^{-b+b^*} (3R^2 + R^3) \leq \frac R2, \label{lst_eq8}
\end{align}

we have that $\Gamma_T$ maps $B_{H}(0,R)$ into itself.\\

Let us prove that it is possible to choose $T\in(0,\frac12]$ in such a way that $\Gamma_T: B_{H}(0,R) \to B_{H}(0,R)$ is a contraction. Let us consider $(u,v)$ and $(w,z)$ in $B_{H}(0,R)$. Then, by using again Lemmas \ref{le2}, \ref{le3}, \ref{ble1}, \ref{ble2}, \ref{ble3} and \ref{ble4}, we have that

\begin{align}
\notag \|\Gamma_T (u,v) - \Gamma_T (w,z) \|_{H}  = & \left\|\eta(\cdot_t)\int_0^{\cdot_t} U (\cdot_t{\scriptstyle -}t') [F_T((u(t'),v(t')) - F_T(w(t'),z(t'))] dt'\right\|_{X^{s+\beta,b}}\\
\notag & + \left\|\eta(\cdot_t)\int_0^{\cdot_t} V (\cdot_t{\scriptstyle -}t') [G_T((u(t'),v(t')) - G_T(w(t'),z(t'))] dt'\right\|_{Y^{s,b}}\\
\notag \leq & CT^{-b+b^*} (\| \alpha uv + \gamma |u|^2 u -\alpha wz - \gamma |w|^2 w \|_{X^{s+\beta,-b}}\\
\notag& + \| \epsilon \partial_x |u|^2 - \partial_x v^2 - \epsilon \partial_x |w|^2 + \partial_x z^2 \|_{Y^{s,-b}})
\end{align}

\begin{align}
\notag\text{}\hspace{4cm}\leq & C T^{-b+b^*} (\|\alpha v(u-w) +\alpha w(v-z)\|_{X^{s+\beta,-b}}\\
\notag &+ \|\gamma |u^2|(u-w) + \gamma w(|u|^2-|w|^2)\|_{X^{s+\beta,-b}}\\
\notag&+ \| \epsilon \partial_x[ |u|^2-|w|^2]\|_{Y^{s,-b}} + \|\partial_x[(z+v)(z-v)]\|_{Y^{s,-b}})\\
\notag \leq & C T^{-b+b^*} ( \|u-w\|_{X^{s+\beta,b}} \|v\|_{Y^{s,b}} + \|w\|_{X^{s+\beta,b}} \|v-z\|_{Y^{s,b}}\\
\notag & +\|u\|^2_{X^{s+\beta,b}} \|u-w\|_{X^{s+\beta,b}} + \|w\|_{X^{s+\beta,b}} \| |u|+|w|\|_{X^{s+\beta,b}} \||u|-|w|\|_{X^{s+\beta,b}}\\
\notag & + \| |u|+|w|\|_{X^{s+\beta,b}} \| |u|-|w| \|_{X^{s+\beta,b}} + \|z+v\|_{Y^{s,b}} \|z-v\|_{Y^{s,b}})\\
\notag \leq & CT^{-b+b^*} (R \|u-w\|_{X^{s+\beta,b}} + R \|v-z\|_{Y^{s,b}}+ R^2 \|u-w\|_{X^{s+\beta,b}}\\
\notag &  + 2R^2 \|u-w\|_{X^{s+\beta,b}} + 2R \|u-w\|_{X^{s+\beta,b}}+2R \|z-v\|_{Y^{s,b}})\\
\notag \leq & C T^{-b+b^*}(6R\|(u,v)-(w,z)\|_{H} + 3R^2 \|(u,v)-(w,z)\|_H).
\end{align}

If we choose $T>0$ such that it satisfies \eqref{lst_eq8}, and $CT^{-b+b^*} (6R+3R^2)<1$ as well, we conclude that $\Gamma_T:B_{H}(0,R) \to B_{H}(0,R)$ is a contraction. Hence there exist $T\in (0,\frac12]$ and a unique couple $(u,v)\in B_H(0,R)$ such that $\Gamma_T(u,v)=(u,v)$.\\

Now we prove that $u\in C(\mathbb R_t; H^{s+\beta}(\mathbb R_x))$ and that $v\in C(\mathbb R_t;H^s(\mathbb R_x))$. Let us note that 
\begin{align}
u(t) = \eta(t) U(t) u_0 - i \eta(t) \int_0^t U(t-t') F_T(u(t'),v(t')) dt'. \label{lst_eq9}
\end{align}

Since $u_0\in H^{s+\beta}(\mathbb R_x)$ and $\{U(t)\}_{t\in\mathbb R}$ is a group of linear operators defined on $H^{s+\beta}(\mathbb R_x)$, it is clear that $\eta(\cdot t) U(\cdot t) u_0\in C(\mathbb R_t;H^{s+\beta}(\mathbb R_x))$. To demonstrate that the second term on the right hand side of \eqref{lst_eq9} is continuous from $\mathbb R_t$, with values in $H^{s+\beta}(\mathbb R_x)$, it is enough to see that this term belongs to $X^{s+\beta,\widetilde b}$ for some $\widetilde b>\frac12$.\\

Let us take $\widetilde b>\frac12$ such that $b^*+\widetilde b\leq 1$. Then
$$-\frac12<-b^*<0\leq \widetilde b\leq (-b^*)+1.$$
By applying Lemma \ref{le2} with $T=1$, and Lemmas \ref{le3}, \ref{ble1} and \ref{ble3}, we conclude that
\begin{align*}
\left\|\eta(\cdot_t)\int_0^{\cdot_t} U (\cdot_t{\scriptstyle -}t')F_T(u(t'),v(t')) dt'\right\|_{X^{s+\beta,\widetilde b}}&\leq C \|\eta \left(\tfrac{\cdot_t}{2T} \right) (\alpha uv +\gamma |u|^2 u) \|_{X^{s+\beta,-b^*}}\\
& \leq C T^{-b+b^*} \| \alpha uv +\gamma |u^2|u\|_{X^{s+\beta,-b}}\\
& \leq C T^{-b+b^*}(\|u\|_{X^{s+\beta,b}} \|v\|_{Y^{s,b}} + \|u\|^3_{X^{s+\beta,b}} ) <\infty.
\end{align*}

In this way, we have shown that $u\in C(\mathbb R_t; H^{s+\beta}(\mathbb R_x))$. Similarly, it can be established that $v\in C(\mathbb R_t; H^{s}(\mathbb R_x))$. This completes the proof of Theorem \ref{localsol}.\qed

\begin{remark} In Theorem \ref{localsol}, the restriction of the fixed point $(u,v)$ of $\Gamma_T$ to $\mathbb R\times[0,T]$ provides a local solution in time to the IVP \eqref{maineq}.
\end{remark}

\section{Regularity gain of the nonlinear part of the solution of IVP \eqref{maineq}. Proof of Theorem \ref{more_reg}}

We will only prove inequality \eqref{more_reg_eq1}, being the proof of \eqref{more_reg_eq2} similar. Let $b^*$ and $\widetilde b$ be such that $b<b^*<\frac12 < \widetilde b$, and $b^*+\widetilde b\leq 1$. By using the immersion $X^{s+\beta+\overline a,\,\widetilde b} \hookrightarrow C_b(\mathbb R_t; H^{s+\beta+\overline a}(\mathbb R_x))$ and taking into account that for $0\leq t'\leq t\leq T\leq\frac12$, $\eta(t)=1$ and
$$F_T(u(t'),v(t'))=[\alpha u(t')v(t')+\gamma |u(t')|^2 u(t')],$$
we have that
\begin{align}
\notag&\left \|\int_0^{\cdot_t} U (\cdot_t{\scriptstyle -}t') [-\alpha u(t')v(t')+\gamma |u(t')|^2 u(t')] dt' \right\|_{C([0,T];H^{s+\beta+\overline a}(\mathbb R_x))}\\
\notag& = \left\| \eta(\cdot_t) \int_0^{\cdot_t} U (\cdot_t{\scriptstyle -}t')F_T(u(t'),v(t')) dt'  \right\|_{C([0,T];H^{s+\beta+\overline a}(\mathbb R_x))}\\
\notag &\leq \left\| \eta(\cdot_t) \int_0^{\cdot_t} U (\cdot_t{\scriptstyle -}t')F_T(u(t'),v(t')) dt'  \right\|_{C_b(\mathbb R_t;H^{s+\beta+\overline a}(\mathbb R_x))}\\
& \leq C\left\| \eta(\cdot_t) \int_0^{\cdot_t} U (\cdot_t{\scriptstyle -}t')F_T(u(t'),v(t')) dt'  \right\|_{X^{s+\beta+\overline a,\,\tilde b}}.\label{PT1_2-1}
\end{align}

Applying Lemma \ref{le2} with $-\frac12<-b^*\leq 0\leq \tilde b\leq -b^*+1$, it follows that
\begin{align}
\left\| \eta(\cdot_t) \int_0^{\cdot_t} U (\cdot_t{\scriptstyle -}t')F_T(u(t'),v(t')) dt'  \right\|_{X^{s+\beta+\overline a,\,\tilde b}} \leq C \|F_T(u(\cdot),u(\cdot))\|_{X^{s+\beta+\overline a,-b^*}}.\label{PT1_2-2}
\end{align}

From Lemma \ref{le3} with $b_1:=-b^*$, $b_2:=-b$, we have that
\begin{align}
\| F_T(u(\cdot),v(\cdot)) \|_{X^{s+\beta+\overline a,-b^*}} \leq CT^{-b+b^*} \|\alpha uv+\gamma |u|^2 u \|_{X^{s+\beta+\overline a,-b}}. \label{PT1_2-3}
\end{align}

Since $\beta+\overline a\leq \frac52 b-\frac58$, by using Lemma \ref{ble3} with $a:=\beta+\overline a$, we obtain
\begin{align}
\|\alpha uv\|_{X^{s+\beta+\overline a,-b}}\leq C\|u\|_{X^{s+\beta,b}} \|v\|_{Y^{s,b}},\label{PT1_2-3.2}
\end{align}

and since $\overline a\leq \frac52 b-\frac58 -\beta\leq 6b-\frac52$ (this inequality holds for $b\geq\frac{11}{28}$, and in particular for $\frac9{20}<b<\frac12$), by Lemma \ref{ble1},
\begin{align}
\|\gamma |u^2|u\|_{X^{s+\beta+\overline a,-b}} \leq C \|u\|^3_{X^{s+\beta,b}}. \label{PT1_2-3.3}
\end{align}

From \eqref{PT1_2-1}, \eqref{PT1_2-2}, \eqref{PT1_2-3}, \eqref{PT1_2-3.2} and \eqref{PT1_2-3.3} we can conclude that
\begin{align*}
\left \|\int_0^{\cdot_t} U (\cdot_t{\scriptstyle -}t') [\alpha u(t')v(t')+\gamma |u(t')|^2 u(t')] dt' \right\|_{C([0,T];H^{s+\beta+\overline a}(\mathbb R_x))}\leq C (\|u\|_{X^{s+\beta,b}}\|v\|_{Y^{s,b}}+\|u\|^3_{X^{s+\beta,b}}),
\end{align*}
which completes the proof of \eqref{more_reg_eq1}.\\

The proof of estimate \eqref{more_reg_eq2} is analogous to that of estimate \eqref{more_reg_eq1}, taking into account that, since $\beta<\frac52b-\frac58=2b-\frac14+(\frac b2-\frac38)<2b-\frac14$, by Lemma \ref{ble4},
\begin{align*}
\|\partial_x |u|^2\|_{Y^{s+a,-b}} \leq C \|u\|^2_{X^{s+\beta,b}},
\end{align*}

and since $0\leq a\leq 5\beta-\frac94<5(\frac52b-\frac58)-\frac94=10b-4+(\frac52 b-\frac{11}8)<10b-4$, by Lemma \ref{ble2},
\begin{align*}
\|\partial_x v^2\|_{Y^{s+a,-b}}\leq C \|v\|^2_{Y^{s,b}}.
\end{align*}

This completes the proof of Theorem \ref{more_reg}.

\section{Dispersive blow-up. Proof of Theorem \ref{blowup}}

We start this section by introducing some necessary definitions. \\

For $0\leq \alpha\leq 1$, we say that $f\in C^{0,\alpha}(\mathbb R)$, if $f:\mathbb R\to\mathbb C$ is continuous and
\begin{align*}
\|f\|_{C^{0,\alpha}(\mathbb R)}:=&\|f\|_{C^{0}(\mathbb R)}+\sup_{x,y\in\mathbb R, x\neq y} \frac{|f(x)-f(y)|}{|x-y|^\alpha}\\
\equiv& \sup_{x\in\mathbb R} |f(x)| +\sup_{x,y\in\mathbb R, x\neq y} \frac{|f(x)-f(y)|}{|x-y|^\alpha}<\infty.
\end{align*}

Moreover, it can be seen that the condition $0\leq \alpha \leq \beta\leq 1$ implies $C^{0,\alpha}\supset C^{0,\beta}$, with a continuous inclusion map. Now, for $k\geq 0$, $k\in\mathbb Z$ and $0\leq\alpha\leq 1$, we define $C^{k,\alpha}(\mathbb R)$ as the subspace of functions $f\in C^k(\mathbb R)$ such that
$$\|f\|_{C^{k,\alpha}(\mathbb R)}:=\sum_{j=0}^k \|f^{(i)}\|_{C^{0,\alpha}(\mathbb R)}<\infty.$$

For $0\leq \alpha\leq \beta\leq 1$ and $k\geq0$, $k\in\mathbb Z$, the following inclusions hold:
$$C^{k+1}(\mathbb R)\subset C^{k,\beta}(\mathbb R)\subset C^{k,\alpha}(\mathbb R)\subset C^k(\mathbb R),$$
with continuous inclusion maps. Additionally, if $s>\frac12$ and $0\leq \alpha\leq \min\{s-\frac12,1\}$, then $H^s(\mathbb R)\subset C^{0,\alpha}(\mathbb R)$ with a continuous inclusion map.\\

Now we construct the initial data that we will use in Theorem \ref{localsol} in order to obtain the dispersive blow-up result.\\

Given $\theta>0$ and $x_0\in\mathbb R$, the function $u_0$ defined by
\begin{align}
u_0(x):=\frac{e^{-\frac i{4\theta}(x-x_0)^2}}{(1+x^2)^{5/4}}\label{u_0}
\end{align}
is such that $u_0\in C^{\infty}(\mathbb R)\cap H^{2^-}(\mathbb R)\cap L^\infty(\mathbb R)$. Besides, if $[U(\cdot_t)u_0](\cdot_x)$ is the solution of the IVP
\begin{align}
\left. \begin{array}{rl}
i\partial_t u+\partial_x^2 u &\hspace{-2mm}=0,\\
u(x,0)&\hspace{-2mm}= u_0(x),
\end{array} \right\}\label{Schr_lin}
\end{align}
then
\begin{itemize}
\item[(i)] for $t\in \mathbb R\setminus \{\frac 1{4\theta}\}$, $[U(\cdot_t)u_0](\cdot_x)\in C^{1,\frac12+\epsilon}(\mathbb R)$, for each $\epsilon\in(0,\frac12)$,
\item[(ii)] $[U(\frac1{4\theta})u_0](\cdot_x)\in C^{1,\frac12+\epsilon}(\mathbb R\setminus\{x_0\})$ for every $\epsilon\in(0,\frac12)$,
\item[(iii)] $[U(\frac1{4\theta})u_0](\cdot_x)\notin C^{1,\frac12+\epsilon}(\mathbb R)$ for every $\epsilon\in(0,\frac12)$.
\end{itemize}

For a proof of the above statements we refer to \cite{LP2019}, \cite{BPSS2014} and \cite{BS2010}.\\

On the other hand, following the ideas presented in \cite{BJM2025} and \cite{LPS2017}, we observe that the function $v_0$ defined by
\begin{align}
v_0(x):=\sum_{j=1}^\infty \alpha_j [V(-\tfrac1{4\theta} j )\phi](x), \label{v_0}
\end{align}

where
$$\phi(x)=e^{-2|x|}\quad \text{and}\quad \sum_{j=1}^\infty \alpha_j e^{4j}<\infty,$$
is such that $v_0\in C^{\infty}(\mathbb R)\cap H^{\frac32^-}(\mathbb R)\cap L^\infty(\mathbb R)$. Furthermore, if $[V(\cdot_t)v_0](\cdot_x)\in C(\mathbb R_t;H^{\frac32^-}(\mathbb R))$ is the solution of the IVP
\begin{align}
\left. \begin{array}{rl}
\partial_t v+\partial_x^5 v &\hspace{-2mm}=0,\\
v(x,0)&\hspace{-2mm}= v_0(x),
\end{array} \right\}\label{KdV5_lin}
\end{align}
then
\begin{itemize}
\item[(i)] for $t\in \mathbb R\setminus \frac 1{4\theta}\mathbb N$, $[V( t)v_0](\cdot_x)\in C^1(\mathbb R)$,
\item[(ii)] $[V(t)v_0](\cdot_x)\in C^1(\mathbb R\setminus\{0\})$ for every $t\in\frac1{4\theta}\mathbb N$,
\item[(iii)] $[V(t)v_0](\cdot_x)\notin C^1(\mathbb R)$ for every $t\in\frac1{4\theta}\mathbb N$.
\end{itemize}

Now we consider  $\delta\in(0,\frac12)$, which will be determined later. We have that $u_0\in H^{2-\delta}$. Let us consider $b$ and $\beta$ as in Theorem \ref{localsol}, i.e., $\frac9{20}<b<\frac12$, $\frac12<\beta<\frac52b-\frac58$, and let us define $s:=2-\delta-\beta$. It is clear that $0\leq s<\frac32$ and that $0\leq s+\beta<2$. Therefore, $u_0\in H^{s+\beta}$ and $v_0\in H^{s}$. Let $(u,v)$ be the local solution in time of the IVP \eqref{maineq} with initial data $(u_0,v_0)$, guaranteed by Theorem \ref{localsol}, and let $[0,T]$ ($T>0$) be the existence interval of that solution. Then, by Theorem \ref{more_reg}, for $t\in[0,T]$,
\begin{align}
u(t)=&U(t)u_0+u_1(t),\label{locsol_u}\\
v(t)=&V(t)v_0+v_1(t),\label{locsol_v}
\end{align}

where $u_1\in C([0,T]; H^{s+\beta+\overline a}(\mathbb R_x))$ and $v_1\in C([0,T];H^{s+a}(\mathbb R_x))$, with $\overline a=\frac52b-\frac58-\beta$ and $a=5\beta-\frac94$.\\

It can be shown that there exist $b\in(\frac9{20},\frac12)$ and $\beta\in(\frac12,\frac52b-\frac58)$ such that $s+\beta+\overline a>2+\frac18-2\delta$ and $s+a>2-2\delta-\frac14$. Choosing $\delta\in(0,\frac12)$ such that $2+\frac18-2\delta>2+\frac1{16}$ and $2-2\delta-\frac14>\frac32+\frac18$, then for $t\in[0,T]$, $u_1(t)\in H^{s_1}$ with $s_1>2+\frac1{16}$ and $v_1(t)\in H^{\tilde s_1}$ with $\tilde s_1>\frac32+\frac18$.\\

Taking into account that for $s>0$ and $0\leq \alpha\leq \min\{s-\frac12,1\}$, $H^s(\mathbb R)\hookrightarrow C^{0,\alpha}(\mathbb R)$, we deduce that $u_1(t)\in C^{0,\alpha}(\mathbb R)$ for $0\leq \alpha\leq 1$ and $u_1'(t)\in C^{0,\alpha}(\mathbb R)$ for $0\leq \alpha\leq\min\{s_1-1-\frac12,1\}$.\\

Since $\min\{s_1-1-\frac12,1\}\geq \min\{1+\frac1{16}-\frac12,1\}=\frac12+\frac1{16}$, we conclude that $u_1(t)\in C^{1,\alpha}(\mathbb R)$ for $0\leq \alpha\leq \frac12+\frac1{16}$ for every $t\in[0,T]$.\\

On the other hand, since $\tilde s_1>\frac32+\frac18$, we have that $v_1(t)\in C^1(\mathbb R)$ for every $t\in[0,T]$. Let us take $\theta>0$ such that $t^*:=\frac1{4\theta}\in[0,T]$. Hence, from \eqref{locsol_u}, taking into account that $U(t^*)u_0\notin C^{1,\frac12+\epsilon}(\mathbb R)$ for $\epsilon\in(0,\frac12)$ and $u_1(t^*)\in C^{1,\alpha}(\mathbb R)$ for $0\leq \alpha\leq \frac12+\frac1{16}$, we conclude that $u(t^*)\notin C^{1,\alpha}(\mathbb R)$ for $0\leq \alpha\leq\frac12+\frac1{16}$.\\

In this way, there exist $\epsilon\in(0,\frac12)$ such that $u(t^*)\notin C^{1,\frac12+\epsilon}$.\\

Finally, from \eqref{locsol_v}, observing that $V(t^*)v_0\notin C^1(\mathbb R)$ and that $v_1(t^*)\in C^1(\mathbb R)$, we conclude that $v(t^*)\notin C^1(\mathbb R)$. This completes the proof of Theorem \ref{blowup}.\qed\\

\textbf{Acknowledgments}\\

This work was partially supported by Universidad Nacional de Colombia, Sede-Medellín– Facultad de Ciencias – Departamento de Matemáticas – Grupo de investigación en Matemáticas de la Universidad Nacional de Colombia Sede Medellín, carrera 65 No. 59A -110, post 50034, Medellín Colombia. Proyecto: Análisis no lineal aplicado a problemas mixtos en ecuaciones diferenciales parciales, código Hermes 60827. Fondo de Investigación de la Facultad de Ciencias empresa 3062.

\end{document}